\documentclass[11pt]{amsart}
\usepackage{amsmath, amssymb} 

	\normalbaselines						  %

\numberwithin{equation}{section}

\newtheorem{thm}{Theorem}[section]

\newtheorem{cor}[thm]{Corollary}
\newtheorem{lm}[thm]{Lemma}
\newtheorem{conj}[thm]{Conjecture}

\theoremstyle{remark}
\newtheorem{rem}[thm]{Remark}

\newtheorem*{rem*}{Remark}

\newenvironment{entry}
{\begin{list}{X}%
  {%
      \setlength{\labelwidth}{55pt}%
      \setlength{\leftmargin}{\labelwidth}
      \addtolength{\leftmargin}{\labelsep}%
      \setlength{\itemsep}{.4pc}
   }%
}%
{\end{list}}

\newcommand{\pbar}{\ensuremath{\bar{\partial}}}

\newcommand{\re}{\operatorname{Re}}

\newcommand{\T}{\mathbb{T}}
\newcommand{\fS}{\mathfrak{S}}
\newcommand{\f}{\varphi}

\newcommand{\e}{\varepsilon}

\newcommand{\D}{\mathbb{D}}
\newcommand{\C}{\mathbb{C}}
\newcommand{\R}{\mathbb{R}}

\newcommand{\db}{\overline\partial}
\newcommand{\p}{\partial}
\newcommand{\pb}{\overline\partial}

\newcommand{\cE}{\mathcal{E}}

\newcommand{\nm}{\,\rule[-.6ex]{.13em}{2.3ex}\,}

\newcommand{\norm}[1]{\ensuremath{\|#1\|}}

\newcommand{\La}{\langle }
\newcommand{\Ra}{\rangle }

\newcommand{\fdot}{\,\cdot\,}

\newcommand{\cH}{\mathcal{H}}

\newcommand{\cP}{\mathcal{P}}
\newcommand{\cC}{\mathcal{C}}

\newcommand{\bP}{\mathbb{P}}

\newcommand\clos{\operatorname{clos}}
\newcommand\spn{\operatorname{span}}
\newcommand\ran{\operatorname{Ran}}

\newcommand{\ti}[1]{_{\scriptstyle \text{\rm #1}}}

\newcommand{\shto}{\raisebox{.3ex}{$\scriptscriptstyle\rightarrow$}\!}

\begin{document}

\title[Ideals: beyond $3/2$]
{The problem of ideals of $H^\infty$: beyond the exponent $3/2$}

\author{Sergei Treil}
\thanks{This paper is based upon work supported by the National Science Foundation under Grant  DMS-0501065.
}
\address{Department of Mathematics, Brown University, 151 Thayer 
Str./Box 1917,      
 Providence, RI  02912, USA }
\email{treil@math.brown.edu}
\urladdr{http://www.math.brown.edu/\~{}treil}

\date{October, 2006}

\begin{abstract} 
The  paper deals  with the problem of ideals of $H^\infty$: describe increasing functions $\f\ge 0$ such that for all bounded analytic functions $f_1, f_2, \ldots, f_n, \tau$ in the unit disc $\D$ the condition 
$$
|\tau(z) | \le \f \left( \left(\sum |f_k(z)|^2\right)^{1/2} \right) \qquad \forall z\in\D
$$
implies that $\tau$ belong to the ideal generated by $f_1, f_2, \ldots, f_n$, i.e.~that there exist bounded analytic functions $g_1, g_2, \ldots, g_n$ such that $\sum_{k=1}^n f_k g_k =\tau$. 

It was proved earlier by the author that the function $\f(s) =s^2$ does not satisfy this condition. The strongest known positive result in this direction  due to J.~Pau states that the function $\f(s) = {s^2}/({(\ln s^{-1})^{3/2} \ln\ln s^{-1} })$ works. However, there was always a suspicion that the critical exponent at $\ln s^{-1}$ is $1$ and not $3/2$. 

This suspicion turned out (at least partially) to be true, $3/2$ indeed is not the critical exponent. The main result of the paper is that one can take for $\f$ any function of form $\f(s) =s^2 \psi (\ln s^{-2})$, where $\psi :\R_+\to \R_+ $ is a bounded non-increasing function satisfying $\int_0^\infty \psi(x) \,dx <\infty$. In particular any of the functions 
$$\f(s) = s^2/( (\ln s^{-2}) (\ln\ln s^{-2}) \ldots (\underbrace{\ln\ln\ldots \ln}_{m \text{ times}} s^{-2})
(\underbrace{\ln\ln\ldots \ln}_{m+1 \text{ times}} s^{-2})^{1+\e}), \qquad \e>0
$$ 
works.

\end{abstract}

\maketitle

\setcounter{tocdepth}{1}
\tableofcontents

\section*{Notation}

\begin{entry}

\item[$\D$] open unit disk in the complex plane $\C$,
$\D:=\{z\in\C\,:\,|z|<1\}$;

\item[$\T$] unit circle, $\T:=\partial \D=\{z\in\C\,:\,|z|=1\}$;
 
\item[$m$] normalized ($m(\T)=1$) Lebesgue measure on the unit circle, $dm = \frac{|dz|}{2\pi}$;

\item[$dA$] area measure on $\C$;

\item[$\p, \pbar$] $\p$ and $\pbar$-operators, 
$\p =\frac12(\frac{\p}{\p x} - i \frac{\p}{\p y})$, 
$\pbar =\frac12(\frac{\p}{\p x} + i \frac{\p}{\p y})$;

\item[$\Delta$] ``normalized'' Laplacian, $\Delta =\p \pbar=
=\frac14(\frac{\p^2}{\p x^2} + \frac{\p^2}{\p x^2})$;

\item[$H^2$, $H^\infty$] Hardy classes of analytic functions, 
$$
H^p := \left\{ f\in L^p(\T) : \hat f (k) := \int_\T f(z) z^{-k}
\frac{|dz|}{2\pi} = 0\ \text{for } k<0\right\}. 
$$
Hardy classes can be identified with spaces of analytic in the unit disk
$\D$ functions: in particular, $H^\infty$ is the space of all bounded
analytic in $\D$ functions; 

\item[$\norm{\cdot}$, $\nm\,.\,\nm$]  norm; since we are dealing with vector- and
operator-valued functions, we will use the symbol $\|\,.\,\|$ (usually with a
subscript) for the norm in a functions space, while $\nm\,.\,\nm$ is used for the
norm in the underlying vector (operator) space.
 Thus for a vector-valued function
$f$ the symbol
$\|f\|_2$ denotes its
$L^2$-norm, but the symbol $\nm f\nm$ stands for the scalar-valued function whose
value at a point $z$ is the norm of the vector $f(z)$; 

\item[$\La\fdot, \fdot\Ra$] inner product in a Hilbert space;
%

\item[$H^2_E$] vector-valued  Hardy class $H^2$ with values in $E$; 

\item[$L^\infty_{\!E\to E_*}$] class of bounded functions on the unit
circle $\T$ whose values are bounded operators from $E$ to $E_*$; $E$ and $E_*$ here are some separable Hilbert spaces. Note, that $E_*$ is not the \emph{dual} of $E$, which we denote by $E^*$. 

\item[$H^\infty_{\!E\to E_*}$] operator Hardy class  of bounded analytic
 functions whose values are bounded  operators from $E$ to $E_*$; 
$$
\|F\|_\infty := \sup_{z\in \D}
\nm F(z)\nm=\underset{\xi\in
\T}{\operatorname{esssup}}\nm F(\xi)\nm;
$$

\end{entry}

Throughout the paper all Hilbert spaces are assumed to be separable. We
always assume that in any Hilbert space an orthonormal basis is fixed, so
any operator $A:E\to E_*$ can be identified with its matrix. Thus besides
the usual involution $A\mapsto A^*$ ($A^*$ is the adjoint of $A$), we 
have two more: $A\mapsto A^T$ (transpose of the matrix) and $A\mapsto
\overline A$ (complex conjugation of the matrix), so $A^* =(\overline
A)^T =\overline{A^T}$. Although everything in the paper can be presented in
invariant, ``coordinate-free'', form use of transposition and complex
conjugation makes the notation easier and more transparent. 

We also assume a ``linear algebra'' notation. A Hilbert space $E$ will be treated as   a space of column-vectors, and its dual $E^*$  (with respect to the \emph{linear} duality, without complex conjugation), will be  identified with the space of row vectors. 

We will also identify a Hilbert space $E$ with the space of operators $\C\to E$; then for a vector $v$ symbol $v^*$ denotes the linear functional $\La\,\cdot\,,v\Ra$, so the inner product $\La u,v \Ra$ can be rewritten as $\La u,v \Ra = v^*u$. Again, this is in complete agreement with the stadard notation accepted in linear algebra. 

Finally, we will use symbol $E_*$ to denote an auxiliary Hilbert space. The reader should not confuse it with $E^*$ which is the dual of $E$ (the space of row vectors). 

\setcounter{section}{-1}

\section{Introduction and main result}

The paper is devoted to the following problem of ideals: given a functions $f_1, f_2, \ldots, f_n, \tau \in H^\infty$ find when $\tau$ belongs to the ideal generated by $f_1, f_2, \ldots, f_n$, i.e.~to find when there exist functions $g_1, g_2, \ldots, g_n\in H^\infty$ solving the Bezout equation
$$
\sum_k f_k g_k=\tau. 
$$
Using vector notation we can define $f(z)=(f_1(z), f_2(z), \ldots, f_n(z))^T \in \C^n =:E$, $g(z)=(g_1(z), g_2(z), \ldots, g_n(z))\in E^*$, so the Bezout equation is rewritten as $g(z)f(z)=\tau(z)$. Clearly, the condition 
$$
|\tau(z)|\le C \nm f(z)\nm \qquad \forall z\in \D
$$
is necessary for the solvability of the Bezout equation. However, this condition, as it was shown long ago by Rao \cite{Rao67}, is clearly not sufficient.

One can ask if a stronger condition 
$$
|\tau(z)|\le C \nm f(z)\nm^p \qquad \forall z\in \D
$$
for $p>1$ is sufficient for the existence of $g\in H^\infty$, $gf =\tau$. it was pretty soon (in the beginning of 80's) understood that this condition is sufficient if $p>2$, and it is not sufficient if $p<2$. The question about $p=2$ (the so-called T.~Wolff's problem) remained open for almost 20 years, until  is was recently settled by the author in \cite{Treil-Ideals-2002}; it was shown there that for $p=2$ the condition is also not sufficient. 

One can consider a bit more general question, namely to ask for which non-decreasing functions $\f$ 
the condition 
$$
|\tau(z)|\le \f(\nm f(z)\nm) \qquad \forall z\in \D
$$
implies the solvability of the equation $gf=\tau$.

In \cite{Lin-ideals-1993} K.~C.~Lin had shown that the answer is affirmative for functions 
$$
\f(s)= \frac{s^2}{(\ln s^{-1})^{3/2 +\e}}, \qquad \e>0. 
$$
Then it was shown by U.~Cegrell \cite{Cegrell1990}  that the answer is affirmative for 
$$
\f(s) = \frac{s^2}{(\ln s^{-1})^{3/2}(\ln\ln s^{-1})^{3/2} \ln\ln\ln s^{-1} }. 
$$
This result was recently improved by J.~Pau, \cite{Pau2005}, who proved that one can take the function 
$$
\f(s) = \frac{s^2}{(\ln s^{-1})^{3/2} \ln\ln s^{-1} } .
$$

The exponent $3/2$ in $(\ln s^{-1})^{3/2}$ always did not look right to the author, the reason being the known estimate in the corona problem. 

Namely, it was known for some time,  that if 
$1\ge\nm f(z)\nm \ge\delta>0$ then one can find $g\in H^\infty_{E^*}$ solving $gf=1$ and such that $\| g\|_\infty \le C \delta^{-2} \ln \delta^{-1}$ (this result first appear in the unpublished preprint by A.~Uchiyama \cite{Uchiyama1981}, see \cite{Nik-shift} for the modern treatment). And it was understood, at least on the level of heuristics, that there is a connection between estimate in the corona problem (equation $gf=1$) and the problem of ideals.%
\footnote{This connection was rigorously used in \cite{Treil-Ideals-2002}, where to prove that the answer is negative for $\f(s)=s^2$ it was shown that the estimate in the corona problem cannot be better than $C \delta^{-2}\ln\ln \delta^{-2}$. Moreover, it can be seen from the proof in \cite{Treil-Ideals-2002} that the answer is negative for any $\f$ such that $\lim_{s\to0+} \f(s) \delta^{-2}\ln\ln \delta^{-2} =0$.} 
So, it was a reasonable to guess that the critical exponent for $\ln s^{-1}$ should be $1$, and not $3/2$.

The theorem below, which is the main result of this paper shows that the critical exponent is indeed at most $1$. 

\begin{thm}
\label{t0.1}
Let $\psi:\R_+\to \R_+$ be a bounded non-increasing function such that \linebreak $\int_0^\infty \psi(x) dx <\infty$, and let the function $\f:[0,1]\to \R_+$ be defined as 
$$
\f(s) := s^2 \psi\left( \ln (s^{-2})\right). 
$$
Then for all $f\in H^\infty_E$, $\|f\|_\infty\le 1$ and $\tau\in H^\infty$ satisfying 
$$
|\tau(z) | \le \f( \nm f(z)\nm \qquad \forall z\in \D, 
$$
there exists $g\in H^\infty_{E^*}$ such that $gf=\tau $. 
\end{thm}
\begin{rem}
It is clear, that only behavior of $\f$ near $0$ (i.e.~the behavior of $\psi $ at $\infty$) is essential for the conclusion of the theorem. This shows, in particular, that the assumption that $\psi$ is bounded is not a restriction; we can always assume that without loss of generality. 
\end{rem}

\begin{rem}
Let $\psi$ behave for large $x$ as 
$$
\psi(x)= \frac1{ x (\ln x) (\ln_2 x) \ldots(\ln_{n-1}x) (\ln_nx)^{1+\alpha}} . 
$$ 
Then it clearly satisfies the assumptions of the theorem, and the function $\f$ behaves near $0$ as
$$
\f(s)= \frac{s^2}{(\ln s^{-2}) (\ln_2 s^{-2}) \ldots(\ln_{n}s^{-2}) (\ln_{n+1} s^{-2})^{1+\alpha}} . 
$$
\end{rem}

\begin{rem}
Note, that is is sufficient to prove Theorem \ref{t0.1} only for inner  functions $f$, i.e.~for $f$ such that $\nm f(z)\nm=1$ a.e~on $\T$. Moreover, it is sufficient to prove it for the so-called co-outer inner functions $f$, i.e.~for such $f$ that the entries of $f$ do not have a common inner divisor. Indeed, let $f=f\ti i f\ti o$ be  the inner-outer factorization (see Section~\ref{s1.1} below) of the function $f$, and let $\theta$ be a common inner divisor of entries $f^k$ of $f$. Define $f_1:= f\ti i/\theta$, $f_2 = \theta f\ti o$ (note that $f_2$ is a scalar-valued function and that $\| f_1(z)\|_\infty \le 1$ because $\|f\|_\infty\le 1$). Then 
$$
\f(\nm f_1 f_2\nm) = \nm f_1\nm^2 |f_2|^2 \psi(\ln(\nm f_1 f_2\nm^{-2})) \le 
\nm f_1\nm^2 |f_2|^2 \psi(\ln(\nm f_1 \nm^{-2})) = |f_2|^2 \f(\nm f_1\nm ), 
$$
so for $\tau$ from the theorem 
$$
|\tau(z) | \le \f(\nm f(z)\nm \le |f_2(z)|^2 \f(\nm f_1(z)\nm). 
$$
Therefore 
$$
|\tau(z)/f_2(z) | \le |f_2(z)| \f(\nm f_1(z)\nm) \le \f(\nm f_1(z)\nm) .
$$
Assuming that Theorem \ref{t0.1} holds for  functions $f$ which are inner and co-outer, we get $g\in H^\infty_{E^*}$ such that $g f_1 = \tau /f_2$ which implies $gf=\tau$. \hfill\qed
\end{rem}

\section{Reduction to the main estimate}

In this section we show that the existence of $g\in H^\infty_{E^*}$ satisfying $gf=\tau$ follows from the boundedness of some bilinear Hankel form. The rest of the paper will be devoted to the proof of this estimate. 

We start by expressing the problem in more geometric terms. But first let us quickly remind the reader some facts about operator- (matrix-) valued analytic functions. 

\subsection{Functions in $H^\infty_{E\to E_*} $ 
and their inner-outer factorization} 
\label{s1.1}
The fact presented here are well known and can be found in monographs \cite{Helson-lect}, \cite{Nagy-Foias-Book-1970}, \cite{Nik-book-v1}. 

Let us recall that a function $F\in H^\infty_{E\shto E_*}$ is called 
\emph{inner} if $F^*(z)F(z)=I$ a.e.~on $\T$ (i.e.~if its boundary values are \emph{isometries}), and it is called \emph{outer} if the set $F H^2_E:=\{Ff: f\in H^2_E\}$ is dense in $H^2_{E_*}$. 

By the famous Beurling--Lax--Phillips theorem any $z$-invariant subspace $\cE\subset H^2$, $z\cE \subset \cE$ can be represented as $\Theta H^2_{E_*}$, where $\Theta \in H^\infty_{E_*\shto E}$ is an inner function, and such inner function is unique up to a constant unitary factor on the right. 

If  for 2 inner functions $\Theta_1 H^2_{E_1} \subset \Theta_2 H^2_{E_2}\subset H^2_E$, then $\Theta_2$ is a left divisor of $\Theta_1$, i.e.~$\Theta_1 = \Theta_2\theta$, where $\theta\in H^\infty_{E_1\shto E_2}$ is an inner function. 

Every operator-valued function $F\in H^\infty_{E\shto E_*}$ admits the inner-outer factorization $F=F\ti i F\ti o$, where $F\ti i\in H^\infty_{E_{**}\shto E_*}$ is an inner function, and $F\ti o\in H^\infty_{E\shto E_{**}}$ is an outer. 
This factorization is unique up to constant unitary factors (on the left of $F\ti o$ and on the right of $F\ti i$). Moreover, 
$$
\clos\{ F H^2_E\} = F\ti i H^2_{E_{**}}
$$

And finally, $F\in H^\infty_{E\shto E_*}$ is called \emph{co-outer} if the function $F^T$ (or, equivalently, the function $z\mapsto (F(\overline z)^*$) is outer. 

All these definitions are much easier to understand for the columns  $f\in H^\infty_E = H^\infty_{\C\shto E}$. Namely, inner functions are exactly the functions $f$ such that $\nm f(z)\nm =1$ a.e.~on $\T$. The outer $f\ti o$ part of a function $f\in H^2_E$  is a scalar-valued function defined by 
$$
f\ti o(z) = \exp\left\{ \int_\T \frac{\xi + z}{\xi - z} \ln \nm f(\xi)\nm  \, dm(\xi) \right\}.
$$
A function $f\in H^2_E$ is co-outer if it cannot be represented as $f= f_1\theta$,  $f_1\in H^2_E$, $\theta$ is a scalar inner function, i.e.~if entries $f^k$ of $f$ do not have a common inner divisor.  

\subsection{Geometric interpretation of the problem}

\begin{lm}
\label{l1.1}
Let $f\in H^\infty_E$ be an inner and co-outer function  (note that in this case  $\nm f(z)\nm \ne 0$ for all $z\in \D$) and let $\tau \in H^\infty$. Then the equation $gf \equiv \tau$ has a solution $g\in H^\infty_{E^*}$ if and only if there exists an operator-valued function $R\in H^\infty_{E\shto E} $,  such that 
$$
R^2 = \tau R, \qquad \ran R(z) = \spn\{f(z)\}, \quad \forall z\in \{z\in \D:\tau(z)\ne 0\} 
. 
$$
\end{lm}

\begin{rem*}
The condition $R^2 = \tau R$ means simply that the values of $\cP(z):=\tau(z)^{-1} R(z)$, defined for $\{z\in \D:\tau(z)\ne 0\}$, are projections (not necessarily orthogonal) onto $\spn\{f(z)\}$. 
\end{rem*}

\begin{proof}[Proof of Lemma \ref{l1.1}]
One direction is trivial.  Namely, let $gf\equiv \tau$. Define $R=fg$. Then $R^2 =fgfg = f\tau g =\tau fg = \tau R$.  

Let us prove the opposite implication. Assuming that there exist an operator-valued function $R$ satisfying assumptions of the lemma, let us construct the solution $g$.  

For any projection $\cP(z)$ onto $\spn\{f(z)\}$ we have $\cP f = f$, so $Rf =\tau f$.

Consider the inner-outer factorization $R=R\ti i R\ti o$, $R\ti i \in H^2_{\C\shto E} = H^2_E$, $R\ti o\in H^\infty_{E\shto \C} = H^\infty_{E^*}$. Since $R f=\tau f$, we have
$$
\tau f H^2 \subset R H^2_E \subset \clos R H^2_E = R\ti i H^2. 
$$
Note, that the inner part of $\tau f$ is $\tau\ti i f$, where $\tau\ti i$ is the inner part of $\tau$. Therefore, comparing inner parts of $\tau f$ and $R$ we get $\tau\ti i f = \theta R\ti i$, where  $\theta $ is some inner function (both $\tau\ti i$ and $\theta$ are scalar-valued). 
We know that $f$ is co-outer, so $\theta$ has to be a divisor of $\tau\ti i$, $\tau\ti i /\theta \in H^\infty$. 

Since $R\ti i = f\tau\ti i /\theta$ and $Rf =\tau f$, 
$$
f (\tau\ti i /\theta) R\ti o f = \tau f. 
$$
Then 
$$
(\tau\ti i /\theta) R\ti o f = \tau,
$$
i.e.~$g=(\tau\ti i /\theta) R\ti o$ is a solution of the Bezout equation $gf=\tau$. 
\end{proof}

\subsection{Reduction to the main estimate}
To construct the function $R$ let us consider the orthogonal projection $\Pi(z)$ onto $\spn\{f(z)\}$. It is easy to see that a function 
$$
R= \tau \Pi + \Pi V (I-\Pi)
$$
where $V$ is an arbitrary operator-valued function,
satisfies $R^2= \tau R$ and $\ran R(z) = \ran \Pi(z) \linebreak= \spn\{f(z)\}$.%
\footnote{Moreover, it is not hard to show, that any such $R$ can be represented in this form. We do not need this for the proof, but it is nice to know that we did not lose anything here}

Direct computations show that $\Pi = f (f^*f)^{-1} f^*$ and $\p \Pi = (I-\Pi) f' (f^*f)^{-1}f^*$, so $\Pi \p \Pi = 0$ and
\begin{equation}
\label{1.2}
(I-\Pi) (\p \Pi) \Pi	= \p\Pi. 
\end{equation}

Let $h_1\in H^2_E$ and let $h_2\in (H^2_E)^\perp$ (i.e.~$\overline h_2 \in zH^2_E$). We get by Green's formula
\begin{align*}
\int\limits_\T \La \tau\Pi h_1, h_2\Ra dm & = \frac2\pi \iint\limits_\D \p\pb \La \tau\Pi h_1, h_2\Ra \log\frac1{|z|}dA(z)
\\ 
& = 
\frac2\pi \iint\limits_\D \p \La ( \tau \pb\Pi) h_1, h_2\Ra \log\frac1{|z|}dA(z)
\end{align*}
Equality \eqref{1.2} can be rewritten as $\Pi \pb \Pi (I-\Pi) = \pb \Pi$, so if we
define $\xi_1:= (I-\Pi)h_1$, $\xi_2:= \Pi h_2$, then 
\begin{align}
\label{1.3}
\int\limits_\T \La \tau\Pi h_1, h_2\Ra dm & = 
\frac2\pi \iint\limits_\D \p \La \Pi(\tau \pb \Pi)(I-\Pi) h_1, h_2\Ra \log\frac1{|z|}dA(z)
\\ 
\notag
& = 
\frac2\pi \iint\limits_\D \p \La (\tau \pb \Pi) \xi_1, \xi_2\Ra \log\frac1{|z|}dA(z) =: L(\xi_1, \xi_2). 
\end{align}
Note, that the bilinear form $L$ is Hankel, meaning that $L(z\xi_1, \xi_2) = L(\xi_1, \overline z \xi_2)$. 

Suppose we are able to show that $L$ is bounded, 
\begin{equation}
\label{1.4}
|L(\xi_1, \xi_2)|\le C \|\xi_1\|_2\|\xi_2\|_2 \qquad \forall \xi_1, \xi_2. 
\end{equation}
Then, applying an appropriate version of the Nehari theorem, see Theorem \ref{t1.2} below, we conclude that there exists a function $V\in L^\infty_{E\shto E}$ such that 
$$
L(\xi_1, \xi_2) = \int\limits_\T \La V \xi_1, \xi_2 \Ra dm \qquad \forall \xi_1, \xi_2. 
$$
Recalling the definition of $L$, see \eqref{1.3}, we conclude that 
$$
\int_\T \La \tau \Pi h_1, h_2\Ra dm = \int_\T \La \Pi V (I-\Pi) h_1, h_2\Ra dm \qquad \forall h_1\in H^2_E, \ \forall h_2 \in (H^2_E)^\perp .
$$
But that exactly means that $R:= \tau \Pi - \Pi V (I-\Pi) \in H^\infty$, so we have constructed the function $R$!

So, to prove the main result it is sufficient to show that the bilinear form $L$ is bounded, i.e.~to prove the estimate \eqref{1.4}

\subsection{A version of the Nehari Theorem}

In this section we present a version of the Nehari Theorem that gives us the existence of a symbol of a the bounded Hankel form $L$, i.e.~ we re going to show that in \eqref{1.4} holds, then  there exists a function $V\in L^\infty_{E\shto E}$, $\|V\|_\infty\le \|L\|$ such that 
$$
L(\xi_1, \xi_2)= \int_\T \La V\xi_1, \xi_2\Ra dm 
$$
for all $\xi_1=(I-\Pi) h_1$, $h_1\in H^2_E$ and for all $\xi_2 = \Pi h_2$, $h_2\in (H^2_E)^\perp$. 

While it is possible to transform the problem so one can apply the classical vectorial Nehari Theorem, a version of the Nehari theorem proved by  S.~Treil and A.~Volberg \cite{TrlVol} that can be applied directly to our situation. 

Let us state this theorem. Let $\mathcal{H}_1$ and $\mathcal{H}_2$ be two separable Hilbert spaces, let $S_1$ be an expanding operator ($\|S_1 x\|\ge\|x\|$ in $\mathcal{H}_1$ and $S_2$ be a contractive operator ($\|S_2\|\le 1$) in $\mathcal{H}_2$ (for our problem at hand we actually will have that $S_1$ and $S_2$ are isometries).  We are given an orthogonal decomposition of $\mathcal{H}_2=\mathcal{H}_2^+\oplus\mathcal{H}_2^-$, and let $S_2\mathcal{H}_2^+\subset\mathcal{H}_2^+$.  Let $\mathbb{P}_+$ and $\mathbb{P}_{-}$ be orthogonal projections in $\mathcal{H}_2$ onto $\mathcal{H}_2^+$ and $\mathcal{H}_2^-$ respectively.  Then a generalized Hankel operator $\Gamma$, $\Gamma:\mathcal{H}_1\to\mathcal{H}_2^-$ is a bounded linear operator satisfying the following relation
\begin{equation}
\label{HankRel}
\Gamma S_1 f = \mathbb{P}_{-} S_2\Gamma f\qquad\forall f\in\mathcal{H}_1.
\end{equation}
A bounded operator $T:\mathcal{H}_1\to\mathcal{H}_2$ satisfying the commutation relation $TS_1=S_2T$ is called a generalized multiplier.  If $T$ is a generalized multiplier, then operator $\Gamma_T:\mathcal{H}_1\to\mathcal{H}_2^-$ defined by 
$$
\Gamma_T f := \bP_- T f , \qquad f \in \mathcal{H}_1
$$
is a generalized Hankel operator. 

\begin{thm}(Treil, Volberg \cite{TrlVol})
\label{t1.2}
Let $S_1$ be an expanding operator and $S_2$ be a contraction.  Given a generalized Hankel operator $\Gamma$ there exists a generalized multiplier $T$ (an operator $T:\mathcal{H}_1\to\mathcal{H}_2$ satisfying $S_2T=TS_1$) such that $\Gamma=\Gamma_T$ and moreover $\norm{\Gamma}=\norm{T}$.
\end{thm}

We apply this theorem to
\begin{align*}
\mathcal{H}_1 & =\clos_{L^2_E}\{(I-\Pi)h:h\in H^{2}_E\},\\
\mathcal{H}_2 &=\clos_{L^2_E}\{\Pi h: h\in L^{2}_E\} \quad \text{and}\\
 \mathcal{H}_2^{-}& =\clos_{L^2_E}\{\Pi h:h\in (H^{2}_E)^\perp\}.
\end{align*}
The operators $S_{1,2}$ are defined by $S_1=M_z\vert_{\mathcal{H}_1}$ and $S_2=M_z\vert_{\mathcal{H}_2}$ where $M_z$ is simply multiplication by the independent variable $z$.  Then clearly, $S_2^* =M_{\bar z}\vert_{\mathcal{H}_2}$ and $S_2^* \cH_2^-\subset \cH_2^-$ so $S_2\cH_2^+\subset \cH_2^+$. 

The bilinear form $L$, defined initially on a dense subset of $\cH_1\times \cH_2^-$ gives rise to a bounded linear operator $\Gamma:\cH_1\to \cH_2^-$, $L(\xi_1, \xi_2) =(\Gamma \xi_1, \xi_2)$. We want to show that $\Gamma$ is a generalized Hankel operator, so Theorem \ref{t1.2} applies. One can see that on the dense set where $L$ is initially defined
\begin{align*}
\La \Gamma S_1\xi_1, \xi_2\Ra &= L (z\xi_1, \xi_2)  = L(\xi_1, \bar z \xi_2) \\ &= \La\Gamma\xi_1, S_2^* \xi_2\Ra 
   = \La S_2 \Gamma \xi_1, \xi_2\Ra = \La\bP_- S_2\Gamma\xi_1, \xi_2\Ra,
\end{align*}
which means that the relation \eqref{HankRel} holds, i.e.,~that $\Gamma $ is a generalized Hankel operator. By Theorem \ref{t1.2} there exists a multiplier $T: \cH_1\to \cH_2$ such that
$$
L(\xi_1, \xi_2)=\La \Gamma \xi_1, \xi_2\Ra = \La T \xi_1, \xi_2 \Ra \qquad \forall \xi_1 \in \cH_1,\ \forall \xi_2 \in \cH_2^-. 
$$
 As one can easily see, in our case any such multiplier is multiplication by a bounded operator-valued function $V$, whose values $V(z)$, $z\in \T$ are bounded operators from $\ran (I-\Pi(z))$ to $\ran \Pi(z)$. Of course, we can always assume the operators $V(z):E\to E$ by defining $V(z)$ to be zero from $(\ran(I-\Pi(z)))^{\perp}$ to $(\ran\Pi(z))^\perp$. \hfill\qed

\section{Uchiyama type lemmas and Carleson measures}

As it is well known, Carleson measures and embedding theorems play important role in the Corona theorem and related problems. In this section we present some  embedding theorems we will need to prove the main estimate.  

\subsection{Carleson measures for $H^2$}
Lemma below is not required for the proof. It is presented only to illustrate the idea of ``correcting factors'', which will be used later in the embedding theorems on hermitian holomorphic vector bundles, see Lemma \ref{l2.4} below. This lemma also might be of an independent interest. 

The reader wanting continue with the proof, can go directly to Lemma \ref{l0.5}

\begin{lm}
\label{l2.1}
Let $M:(-\infty, 0]\to [0,1]$ be a $C^2$ non-decreasing function such that 
\begin{equation}
\label{0.1}
\left(  \begin{array}{cc}
M & M' \\ 
M' & M'' \\ 
\end{array} \right) \ge 0, \qquad \forall x \in (-\infty,0], 
\end{equation}
 Then for any subharmonic  $u\in C^2(\D)\cap C(\clos\D)$, $u\le 0$ and for any $h\in H^2_E$
\begin{align*}
\frac{2}{\pi}\int\limits_\D M'(u(z)) \Delta u (z) \nm h(z)\nm^2 \log\frac{1}{|z|} \,dxdy 
& \le \int\limits_\T \nm h(z)\nm^2 M(u(z)) \, dm(z) \\
&\le \int\limits_\T \nm h(z)\nm^2  dm(z). 
\end{align*}
i.e.~the measure $\frac{2}{\pi}M'(u(z))\Delta u(z)\log\frac{1}{|z|}\,dxdy$ is Carleson with the norm of the embedding operator at most $1$.
\end{lm}

\begin{proof}
Let us first assume that in addition $h$ is $C^2$-smooth in the closed disc  $\clos\D$. Then by Green's formula
\begin{multline}
\label{2.1}
	\int\limits_\T \nm M(u) h\nm^2dm - M(u(0)) \nm h(0)\nm^2  =
	\frac2\pi \iint\limits_\D \Delta \left( M( u(z)) \nm h(z)\nm^2 \right) \log\frac1{|z|}\,dxdy .
\end{multline}
Computing the Laplacian we get 
\begin{multline*}
	\Delta \left( M( u(z)) \nm h(z)\nm^2 \right)  = M'(u(z)) \Delta u \nm h\nm^2 \\
	+
	 M(u) \nm h'\nm^2 + M''(u) |\p u|^2 \nm h\nm^2 + 2\re M'(u) \db u (h',h) .
\end{multline*}
The assumption $\left(\begin{array}{cc}
M & M' \\ 
M' & M'' \\ 
\end{array}\right)\ge 0$ implies that 
$$
 M(u) \nm h'\nm^2 + M''(u) |\p u|^2 \nm h\nm^2 + 2\re M'(u) \db u (h',h) \ge 0, 
$$ 
so we get from \eqref{2.1}
\begin{align*}
\frac{2}{\pi}\int\limits_\D M'(u(z)) \Delta u (z) \nm h(z)\nm^2 \log\frac{1}{|z|} \,dxdy 
& \le \int\limits_\T \nm h(z)\nm^2 M(u(z)) \, dm(z) 
\end{align*}
which proves the lemma assuming additional smoothness of $h$. 

Using standard reasoning, i.e.~considering discs of radius $r<1$ and taking limit as $r\to1-$ we get the lemma without assuming additional smoothness of $h$. 
\end{proof}

The following lemma was presented to the author by F.~Nazarov.

\begin{lm}[F.~Nazarov, personal communication]
\label{l0.5}
Let $\psi $ be as in Theorem \ref{t0.1}, i.e.~let $\psi$ be a bounded non-increasing function, such that $\int_0^\infty \psi(x)\,dx<\infty$. Then there exists an increasing $C^2$-function $M:(-\infty, 0]\to [0,1]$ satisfying the assumptions of Lemma \ref{l1.1} and such that  
$$
\psi(x) \le C M'(-x), \qquad \forall x\in \R_+ ,
$$
for some $C<\infty$. 
\end{lm}

In other words, Lemma \ref{l0.5} states that one can always replace an arbitrary function $\psi$ satisfying the assumptions of Theorem \ref{t0.1}  by the derivative $M'(-x)$ of a function  $M$ satisfying the assumption of Lemma \ref{l2.1}. 

\begin{proof}[Proof of Lemma \ref{l0.5}]
One can estimate 
\begin{align*}
\psi(x) & \le \psi(0) \mathbf 1_{[0,1]} + \int_1^\infty  \frac1r \mathbf 1_{[0, r]} (x) d(-\psi(r))  
\\
& \le e \psi(0) e^{-x} + e \int_1^\infty \frac1r e^{-x/r} d(-\psi(r)) =: K(x);
\end{align*}
here we used the trivial estimate $\mathbf 1_{[0,r]} (x) \le e e^{-x/r}$. Clearly $K(x) \le N'(-x)$, where 
$$
N (x) = e \psi(0) e^x +  e\int_1^\infty e^{x/r} d(-\psi(r)), \qquad x\in(-\infty, 0].
$$

Functions $e^{x/r}$ obviously satisfy \eqref{0.1}, and this condition \eqref{0.1} is preserved under convex combinations. Therefore, trivially, the function $N$ satisfies \eqref{0.1}. 

Function $N$ on $(-\infty, 0]$ is trivially bounded, $0\le N(x) \le C$, so the function $M=C^{-1}N$ satisfies the conclusion of the lemma. 
\end{proof}

The next lemma is not needed for the proof, and presented only as an illustration of the application of Lemma \ref{l2.1}. 

\begin{lm}
\label{l2.3}
Let $f\in H^\infty_{E}$, $\|f\|_\infty\le 1$, and let $\psi$ be a function satisfying the assumption of Theorem \ref{t0.1},  i.e.~$\psi$ be a bounded non-increasing function such that $\int_0^\infty \psi(x) dx<\infty$.  Then the measure 
$$
\frac{\nm f\nm^2 \nm f'\nm^2 -|(f,f')|^2}{\nm f\nm^4 } \psi(\ln \nm f(z)\nm^{-2})\, \log\frac1{|z|}\, dA(z)
$$
is Carleson. 
\end{lm}
\begin{proof}
Direct computations show that 
$$
\Delta \ln(\nm f(z)\nm^2)= \frac{\nm f\nm^2 \nm f'\nm^2 -|(f,f')|^2}{\nm f\nm^4 }. 
$$
By Lemma \ref{l0.5} there exists a function $M$ satisfying the assumptions of Lemma \ref{l2.1} and such that
$$
\psi (-x) \le M'(x) \qquad \forall x\in (-\infty, 0].
$$
Applying Lemma \ref{l2.1} with this $M$ and $u(z) = \ln\nm f(z)\nm^2$ we get the conclusion of Lemma \ref{l2.3}. 
\end{proof}

\subsection{Embedding on holomorphic vector bundles}
\label{s2.2}

Let $E(z)$, $z\in D$ be an analytically varying family of subspaces of a Hilbert space $H$ (i.e. a subbundle of a trivial bundle), and let $\Pi(z)$ be an orthogonal projection onto $E(z)$. In other words, $E(z)$ \emph{locally} can be represented as $E(z) =\ran F(z)$, where $F$ is an analytic  function whose values are operators $E_*\to E$ and such that $F^*F\ge \delta^2 I$ (but we do not assume a uniform  estimate for all $z\in \D$) 

The projection $\Pi(z)$ can be written as $\Pi = F(F^*F)^{-1}F^*$, so $\p \Pi = (I-\Pi) F'(F^*F)^{-1}F^*$. From this formula  it is easy to see that  
the function $\Pi$ satisfied the equation $\Pi \p\Pi =0$.  It is not hard to show that the identity $\Pi\p\Pi=0$  can be used as an equivalent definition of an analytic family of subspaces. However, in what follows we do not need the equivalence. Formally we will only assume that $\Pi$ is a $C^2$-smooth function whose values are orthogonal projections satisfying $\Pi\p\Pi=0$.

To prove that we will need the following lemma. It was proved in \cite{TrWick}, and we present the proof here  only for the convenience of the reader. 

\begin{lm}
\label{PdP}
Let $\Pi$ be a $\cC^2$ smooth function (of one complex variable) whose values are orthogonal projections in a Hilbert space. Assume that $\Pi\p\Pi = 0$. Then  
\begin{align*}
&(\p \Pi) (I-\Pi) =0, \qquad\p \Pi = (\p\Pi) \Pi = (I-\Pi)\p\Pi\qquad\text{and}\\
&\Delta \Pi := \p\pbar \Pi = (\p\Pi)(\p\Pi)^* - (\p\Pi)^*(\p\Pi)
\end{align*}
\end{lm}

\begin{rem}
\label{r2.5}
Since $\db\Pi = (\p\Pi)^*$, by taking conjugates we get the following identities for $\db\Pi$:
$$
(\db\Pi)\Pi = (I-\Pi) \db \Pi =0, \qquad \db\Pi =\Pi \db\Pi = (\db\Pi)(I-\Pi). 
$$
\end{rem}

\begin{proof}[Proof of Lemma \ref{PdP}]
Using  the identity $\Pi=\Pi^2$ we get
$$
\p\Pi = \p \Pi^2 =(\p\Pi)\Pi+\Pi\p\Pi = (\p\Pi)\Pi,
$$
because $\Pi\p \Pi=0$. Thus we have proved that $(\p\Pi)\Pi =\p \Pi$. 
The identity $(\p\Pi)(I-\Pi)=0$ follows immediately because 
$$
(\p\Pi)(I-\Pi) = \p\Pi - (\p\Pi)\Pi =\p\Pi -\p\Pi = 0. 
$$

The identity $(I-\Pi)\p\Pi = \p\Pi$ is an immediate corollary of the hypothesis $\Pi\p\Pi =0$:
$$
(I-\Pi)\p\Pi = \p\Pi -\Pi\p\Pi = \p\Pi. 
$$

Let us now prove formula for $\Delta\Pi$. The taking the adjoints of both sides of identity $\p\Pi = (\p\Pi)\Pi$ and using the fact that $(\p\Pi)^*=\overline\p\Pi$ we get 
$$
\overline\p \Pi = (\p\Pi)^* = \Pi (\p\Pi)^* = \Pi \pb \Pi,
$$
so $\pb \Pi = \Pi \pb \Pi$. Applying $\p$ to both sides of this identity we get
\begin{eqnarray*}
\p\pbar\Pi & = & \p\Pi\pbar\Pi+\Pi\p\pbar\Pi\\
 & = & \p\Pi\pbar\Pi-\pbar\Pi\p\Pi+\pbar\Pi\p\Pi+\Pi\p\pbar\Pi\\
 & =& \p\Pi\pbar\Pi-\pbar\Pi\p\Pi+\pbar\left(\Pi\p\Pi\right).
\end{eqnarray*}
Using the hypothesis that $\Pi\p\Pi=0$ and the fact that $\pb\Pi= (\p\Pi)^*$, see above Remark \ref{2.5}, we get the final identity.
\end{proof}

We are interested in embedding theorem for functions of form $\xi(z) = \Pi(z) h(z)$, where $h\in H^2_-(E)$. 

Note, that such functions are dense in the set of all antiholomorphic (with respect to covariant derivative) sections on the hermitian (holomorphic) vector bundle. We do not use directly this fact in the proof, but this remark might help a reader with the background in geometry to understand better what is going on. 

\begin{lm}
\label{l2.4}
Let a function $M$ satisfies the assumptions of Lemma \ref{l2.1}, and let $u\le 0$ be a subharmonic function satisfying $\Delta u(z) \ge \nm \p \Pi(z) \nm^2$. Then for all $\xi$ of form $\xi(z) = \Pi(z) h(z)$,  $\overline h\in H^2(E)$  the embedding 
\begin{align*}
\frac{2}{\pi}\int\limits_\D e^{u(z)} M'(u(z)) \nm\p \Pi(z)\nm^2 \nm \xi(z)\nm^2 \log\frac{1}{|z|} \,dxdy & \le \int\limits_\T e^{u(z)} M(u(z)) \nm \xi(z)\nm^2  \, dm(z) \\
&\le \int\limits_\T \nm \xi(z)\nm^2  dm(z). 
\end{align*}
holds. To avoid discussion about boundary values we assume here that $\Pi$ and $h$ are continuous up to the boundary. 
\end{lm} 

\begin{proof}
Applying Green's formula to the integral 
$
\int_\T B(u) \nm \xi(z)\nm^2 dm
$ where $B(u) := e^u M(u)$ we get 
\begin{equation}
\label{2.2}
\int_\T B(u) \nm \xi(z)\nm^2 dm -B(u(0))\nm \xi(0)\nm^2 = \iint_\D \Delta (B \nm \xi\nm^2)\log \frac1{|z|}\,dxdy. 
\end{equation}

Let us first compute the Laplacian. We get
$$
\p (B\nm \xi\nm^2) = B' \p u \nm \xi\nm^2 + B \La \p\xi , \xi\Ra + B \La\xi, \pb\xi\Ra . 
$$
Note, that the identity $(\p\Pi )\Pi = \p\Pi$ and the fact that $\p h=0$ imply that 
$$
\p\xi =\p (\Pi h) = (\p\Pi) h = (\p\Pi) \Pi h = (\p\Pi) \xi. 
$$ 
Using this and  the identity $\Pi\p\Pi=0$  we get 
$$
\La \p\xi, \xi\Ra = \La \p\xi,  \Pi \xi \Ra =  \La ( \p\Pi) \xi,  \Pi \xi \Ra =0, 
$$
so 
$$
\p (B\nm \xi\nm^2) = B' \p u \nm \xi\nm^2  + B \La\xi, \pb\xi\Ra
$$
Taking $\pb$ derivative of this equation and using the identity $\La \xi , \p\xi \Ra=0$ we get 
$$
\Delta (B\nm \xi\nm^2) = B' \Delta u \nm \xi\nm^2 +
\left[ B'' |\p u|^2 \nm\xi\nm^2 + 2 \re( B'\p u \La \pb \xi, \xi \Ra )+ B \nm \pb \xi \nm^2 \right] + B\La \xi , \Delta \xi\Ra
$$
To handle the term $\La \xi , \Delta \xi\Ra$ we take the $\p$ derivative of the equation $\La \xi , \p\xi \Ra = 0$ to get 
$$
\La \p \xi, \p\xi \Ra + \La \xi , \Delta \xi\Ra, 
$$
so $\La \xi , \Delta \xi\Ra = -\nm \p \xi\nm^2 = - \nm (\p\Pi) \xi\nm^2$.  Thus we can rewrite the Laplacian as 
\begin{align}
\label{2.3}
\Delta (B\nm \xi\nm^2) = B' \Delta u \nm \xi\nm^2  &+
\left[ B'' |\p u|^2 \nm\xi\nm^2 + 2 \re( B'\p u \La \pb \xi, \xi \Ra )  + B \nm \pb \xi \nm^2 \right]  
\\ \notag
& - B \nm (\p\Pi) \xi\nm^2. 
\end{align}
The expression in brackets is just the quadratic form 
$$
\left\La 
\left( \begin{array}{cc} B & B' \\ B' & B'' \end{array} \right) 
\left( \begin{array}{c} \pb \xi\\ (\pb u) \xi \end{array}\right), 
\left( \begin{array}{c} \pb \xi\\ (\pb u) \xi \end{array}\right)
\right\Ra, 
$$
so it is non-negative if  $\left( \begin{array}{cc} B & B' \\ B' & B'' \end{array} \right) \ge 0$. Recall that $B(u)= e^u M(u)$, so $B' = e^u (M+M')$, $B'' = e^u (M + 2 M' + M'')$. Note that $M, M', M'' \ge 0$, so to show positive semi-definiteness we only need to check the determinant, which is 
$$
e^{2u} \left\{ M(M+2M'+M'') - (M+M')^2\right\} = e^{2u} \left\{  MM'' - (M')^2\right\}\ge 0
$$
The last inequality, together with $M''\ge 0$, follows from the fact that $\left( \begin{array}{cc} M& M'\\ M' & M'' \end{array}\right) \ge 0$.

We can also use the inequality $\Delta u \ge \nm\p \Pi\nm^2$ to estimate the difference 
$$
B' \Delta u \nm \xi\nm^2 - B\nm (\p \Pi) \xi\nm^2 \ge  e^u \left\{ (M+M') \Delta u \nm\xi\nm^2 - M \nm \p\Pi\nm^2 \nm \xi \nm^2 \right\} \ge e^u M'(u) \Delta u \nm\xi\nm^2. 
$$
Therefore we will get from \eqref{2.2}
$$
\int\limits_\T B(u) \nm \xi(z)\nm^2 dm \ge \iint\limits_\D e^u M'(u) \Delta u \nm\xi\nm^2 \log\frac1{|z|} dxdy \ge \iint\limits_\D e^u M'(u) \nm \p\Pi\nm^2 \nm\xi\nm^2 \log\frac1{|z|} dxdy, 
$$
whicn proves the lemma. 
\end{proof}

\begin{lm}
\label{l2.5}
Let $\Pi$, $M$, $u$ be as in Lemma \ref{l2.4}, and let $g\in H^\infty$, $\|g\|_\infty\le 1$  be continuous up to the boundary and has finitely many zeroes in $\D$.  Assume that 
$$
|g(z) | \le e^{u(z)} M'(u(z)). 
$$

Then for any $\xi$ of form $\xi = \Pi h$, $\overline h \in H^2(E)$  (we again assume that $\Pi$ and $h$ are continuous up to the boundary) we get the embedding
\begin{align*}
\frac{2}{\pi}\int\limits_\D  \nm\pb( \overline{g^{1/2} }\xi )\nm^2  \log\frac{1}{|z|} \,dxdy &\le 2\int\limits_\T \nm \xi(z)\nm^2  dm(z). 
\end{align*}
\end{lm}

\begin{rem*}
Let us say few words about interpretation of the integral in the left side. If $f(z)\ne 0$ we can say that in a neighborhood of $z$ $g^{1/2}$ is a branch of a square root. It is easy to see that the expression $\nm \pb (\overline g^{1/2} \xi)\nm$ does not depend on the choice of the branch. And of course, in the integral  we ignore the points where $f(z) =0$. 

Another, more ``high brow'' explanation is that the function $g^{1/2}$ is defined on its Riemann surface. Then the function $ \nm \pb (\overline g^{1/2} \xi)\nm$ is defined on this Riemann surface, and it can be pushed back to a single-valued function on the disc $\D$. 
\end{rem*}

\begin{proof}[Proof of Lemma \ref{l2.5}]
Let us define $\eta:= \overline{g^{1/2} }\xi$, and apply the Green's formula:%
\footnote{There is a delicate moment here to justify the formula, because the Laplacian $\Delta(\nm \eta\nm^2)$ is not defined at zeroes of $g$. So one essentially needs to repeat the proofs of the Green's formula. Namely, one needs to apply the identity $\iint_G (U \Delta V - V\Delta U) dA = \int_{\p G} (U \frac{\p V}{\p n} - V \frac{\p U}{\p n} )dm$ with $U=\log\frac1{|z|}$ and $V=\nm \eta\nm^2$. The domain $G$ here is the unit disc $\D$  without small discs $D_j$ around zeroes of $g$ and a small disc $D_0$ around the origin. 
Since the gradient $\nabla V$ is bounded around zeroes of $g$ (although it is not defined at zeroes of $g$), the integrals $\int_{\p D_j}\ldots$ tend to $0$ as we shrink the radii of $D_j$. The integral $\int_{\p D_0} \ldots$ tends to $-V(0)$, so shrinking the radii and taking th limit we get the equality  that looks exactly like the classical Green's formula. 
}
\begin{equation}
\label{2.4}
\int\limits_\T \nm \eta \nm^2 dm - \nm \eta(0)\nm^2 = \frac{2}{\pi} \iint\limits_\D \Delta (\nm \eta\nm^2) \log\frac1{|z|} dA(z)
\end{equation}
Computing the Laplacian (cf \eqref{2.3} with $B\equiv 1$) we get
$$
 \Delta (\nm \eta\nm^2) = \nm \pb \eta\nm^2 - \nm (\p\Pi) \eta \nm^2 .
$$
Substituting it to \eqref{2.4} we get 
$$
 \frac{2}{\pi} \iint\limits_\D \nm\pb \eta\nm^2 \log\frac1{|z|} dA(z) \le \int\limits_\T \nm \eta \nm^2 dm +
  \iint\limits_\D \nm (\p\Pi) \eta \nm^2  \log\frac1{|z|} dA(z) . 
$$
Noticing that 
$$
\nm (\p\Pi) \eta \nm^2 \le \nm \p\Pi\nm^2 |g| \,\nm \xi \nm^2 \le e^u M'(u) \nm \p\Pi\nm^2  \nm \xi \nm^2,
$$
and applying Lemma \ref{l2.4}, we can estimate the second term as
$$
  \frac{2}{\pi} \iint\limits_\D \nm (\p\Pi) \eta \nm^2  \log\frac1{|z|} dA(z)\le \int\limits_\T \nm \xi(z)\nm^2  dm(z)
$$
which proves the lemma (because $\nm \eta\nm \le\nm \xi\nm$). 
\end{proof}

\begin{cor}
\label{c2.6}
Let $f\in H^\infty_E$, $\|f\|_\infty\le1$, $\nm f(z)\nm \ne 0$ for all $z\in\D$, and let $\Pi(z)$ be the orthogonal projection onto $\spn\{f(z)\}$.

Let $\psi: \R_+\to \R_+$ be a bounded  non-increasing function satisfying $\int_0^\infty \psi(x)\,dx<\infty$. Define  $\f:[0,1]\to\R_+$ as 
$$
\f(s) = s^2 \psi(\ln s^{-2}), 
$$ 
and assume that  $\tau \in H^\infty$ satisfies $|\tau(z)| \le \f(\nm  f (z)\nm )$. 

Then for any $\xi$ of form $\xi = \Pi h$, $\overline h \in H^2_E$
\begin{align}
\label{2.5}
\frac2\pi\iint\limits_\D |\tau| \nm \p\Pi\nm^2 \nm \xi\nm^2 \log\frac1{|z|}\,dA(z)  & \le \int\limits_\T \nm \xi\nm^2 dm, 
\\
\label{2.6}
\frac2\pi\iint\limits_\D \nm \pb (\overline{\tau^{1/2}} \xi) \nm^2 \log\frac1{|z|}\,dA(z)   & \le 2  
\int\limits_\T \nm \xi\nm^2 dm. 
\end{align}

Moreover,  for all $\zeta$ of form $\zeta = (I-\Pi) h$, $h\in H^2_E$
\begin{align}
\label{2.7}
\frac2\pi\iint\limits_\D |\tau| \nm \p\Pi\nm^2 \nm \zeta \nm^2 \log\frac1{|z|}\,dA(z)  & \le \int\limits_\T \nm \zeta\nm^2 dm, 
\\
\label{2.8}
\frac2\pi\iint\limits_\D \nm \p ({\tau^{1/2}} \zeta ) \nm^2 \log\frac1{|z|}\,dA(z)   & \le 2  
\int\limits_\T \nm \zeta\nm^2 dm. 
\end{align}

Here again, to avoid complications with the existence of the boundary values, we  assume that the functions $\Pi$ and $h$ are continuous up to the boundary. 
\end{cor}

\begin{proof}
Direct computations show that 
$$
\nm \p \Pi\nm^2 = \frac{\nm f \nm^2 \nm f'\nm^2 - |\La f', f \Ra|^2}{\nm f \nm^4}. 
$$
Probably the easiest way to see that is to treat everything in non-commutative settings, i.e. consider the case where $f$ is an arbitrary operator-valued function. Then $\Pi = f (f^*f)^{-1} f^*$ and 
$$
\p \Pi = (I-\Pi) f' (f^*f)^{-1} f^*. 
$$ 
The operator $\p\Pi(z)$ is in our case a rank one operator, so its operator norm coincides with its Hilbert--Schmidt norm. The latter can be computed as 
$$
\nm \p\Pi\nm^2_{\fS_2} = \nm f' (f^*f)^{-1} f^* \nm^2_{\fS_2} - 
\nm \Pi f' (f^*f)^{-1} f^* \nm^2_{\fS_2} 
= \frac{\nm f'\nm^2}{\nm f\nm^2} - \frac{|\La f', f\Ra|^2}{\nm f\nm^4}. 
$$

On the other hand, for $u(z)= \log \nm f(z)\nm^2$
$$
\Delta u(z) = \frac{\nm f \nm^2 \nm f'\nm^2 - |\La f', f \Ra|^2}{\nm f \nm^4} = \nm \p \Pi\nm^2  .
$$
This is not a coincidence, because in our case $\nm\p \Pi(z) \nm^2$ is up to the sign ``$-$'' the curvature of the corresponding holomorphic Hermition vector bundle. And the same curvature can be computed by taking the Laplacian of $\log \nm f\nm^{-2}$. 

As it was shown in Lemma \ref{l0.5}, 
$$
\psi(x) \le C M'(-x), \qquad x\in \R_+, 
$$
where $M$ is some function 
satisfying the assumptions of Lemma \ref{l2.1} (and so of Lemmas \ref{l2.4}, \ref{l2.5}).

Then, for  $u=\log \nm f \nm^2$ we get 
\begin{align*}
|\tau| \le \f(\nm f\nm)  & = \nm f\nm^2 \psi(\ln \nm f\nm^{-2}) \\
& = e^u \psi (-u) \\
& \le e^u M'(u), 
\end{align*}
so we are now in position to apply  Lemmas \ref{l2.4}, \ref{l2.5}. Namely,
$$
\frac2\pi\iint_\D |\tau|  \nm \p\Pi\nm^2 \nm \xi\nm^2 \log\frac1{|z|}\,dA(z) \le 
\frac2\pi\iint_\D e^u M'(u) \nm \p\Pi\nm^2 \nm \xi\nm^2 \log\frac1{|z|}\,dA(z) 
 \le \int\limits_\T \nm \xi\nm^2 dm
$$
where the last inequality follows from Lemma \ref{l2.4}. The estimate \eqref{2.6} is obtained by immediate application of Lemma \ref{l2.5}. 

To prove the other 2 inequalities, define $Q=I-\Pi$. Clearly, $Q$ is projection onto antianalytic family of subspaces, namely it satisfies $Q\pb Q=0$. Then $Q(\overline z)$ is a projection onto analytic family of subspaces. Note also that $\p Q = -\p \Pi$, so 
$$
\nm \p Q \nm = \nm \pb Q\nm = \nm \p \Pi\nm. 
$$ 
Therefore,  change of variables $z\mapsto \overline z$ reduces \eqref{2.7}, \eqref{2.8} to the estimates \eqref{2.5}, \eqref{2.6}, which we already proved. 
\end{proof}

\section{Main estimates}

To complete the proof of Theorem, we need to show that the bilinear form is bounded, i.e~to prove the estimate \eqref{1.4}. 

We want to estimate 
$$
L(\xi_1, \xi_2) = \frac2\pi \iint\limits_\D \p \La \tau (\pb \Pi) \xi_1, \xi_2 \Ra\, \log\frac1{|z|} dA(z) 
$$
where $\xi_1 = (I-\Pi) h_1$, $h_1\in H^2_E$, $\xi_2 = \Pi h_2$, $\overline h_2 \in zH^2_E$.

Define $\eta_1 := \tau^{1/2} \xi_1$, $\eta_2:= \overline{\tau^{1/2}}\, \xi_2$, both functions are defined on the Riemann surface of $\tau^{1/2}$. Note, that if $z^\sharp$ is a point on the Riemann surface corresponding to a point $z\in \D$, then 
$$
\La (\pb\Pi(z)) \eta_1(z^\sharp), \eta_2(z^\sharp)\Ra  = \La \tau(z) (\pb\Pi(z)) \xi_1(z) , \xi_2(z)\Ra. 
$$
In particular, this expression does not depend on the choice of $z^\sharp$ corresponding to $z\in \D$, i.e.~on the choice of the branch of $\tau^{1/2}$. 

Let us rewrite $L(\xi_1, \xi_2)$ in terms of $\eta_{1}$, $\eta_2$. To shorten the notation let us denote the measure $\frac2\pi \log\frac1{|z|} dA(z)$ by $\mu$. Then 
\begin{eqnarray*}
L(\xi_1,\xi_2) & = & \int _{\D}\p\La ( \db\Pi) \eta_1,\eta_2\Ra d\mu\\
 & = & \int_{\D}\La\Delta \Pi\eta_1,\eta_2\Ra d\mu+\int_{\D}\La(\db\Pi)\p\eta _1,\eta_2\Ra d\mu+\int_{\D}\La (\db\Pi) \eta_1,\db\eta_2\Ra d\mu\\
 & := & \textnormal{I} + \textnormal{II} + \textnormal{III}.
\end{eqnarray*}

We claim that the first integral disappears, $I=0$.

Lemma \ref{PdP} gives us
$$
\p\db\Pi = \p\Pi\db\Pi-\db\Pi\p\Pi. 
$$
Then
\begin{eqnarray*}
\textnormal{I} & = & \int_{\D}\La\p\db\Pi\eta_1,\eta_2\Ra d\mu\\
 & = & \int_{\D}\La\p\Pi\db\Pi\eta_1,\eta_2\Ra d\mu-\int_{\D}\La\db\Pi\p\Pi\eta_1,\eta_2\Ra d\mu\\
 & = & \int_{\D}\La\db\Pi\eta_1,\db\Pi\eta_2\Ra d\mu-\int_{\D}\La\p\Pi\eta_1,\p\Pi\eta_2\Ra d\mu=0.
\end{eqnarray*}
Using Lemma \ref{PdP} we have that $\p\Pi\eta_1=\p\Pi (I-\Pi) \eta_1 =0$, $\db\Pi\eta_2= \pb \Pi \Pi \eta_2 =0$, and $(\p\Pi)^*=\db\Pi$, so the both summands disappear and $\mathrm{I}=0$.

Let us estimate $\mathrm{II}$:
$$
|\mathrm{II}| \le \int_\D \nm \p \Pi\nm \nm \eta_2 \nm \nm \p \eta_1\nm \,d\mu
\le \left(\int_\D \nm \p \Pi\nm^2 \nm \eta_2 \nm^2 \,d\mu \right)^{1/2}
\left(\int_\D \nm \p \eta_1\nm^2 \,d\mu  \right)^{1/2}. 
$$
Using \eqref{2.7} from Corollary \ref{c2.6} we get
$$
\int_\D \nm \p \Pi\nm^2 \nm \eta_2 \nm^2 \,d\mu = \int_\D \nm \p \Pi\nm^2 |\tau| \cdot\nm \xi_2 \nm^2 \,d\mu \le \int_\T \nm \xi_2\nm^2 dm, 
$$
and by \eqref{2.6} 
$$
\int_\D \nm \p \eta_1\nm^2 \,d\mu  \le 2\int_\T \nm \xi_1\nm^2 \,dm .
$$
The integral $\mathrm{III}$ is estimated similarly, one needs to use \eqref{2.5} and \eqref{2.8} in that case. \hfill\qed

\section{Concluding remarks} 

The author is grateful to F.~Nazarov for helpful discussions, especially for introducing Lemma \ref{l0.5}, which has allowed greatly simplify the statement of the main result. The author is not sure whether this lemma was known before, but the statement and the proof presented in the text belong to F.~Nazarov.  

There were 2 crucial new ideas in this paper that had allowed us to get a better result. 
The first one is a more careful estimate of the Carleson measure. Namely, it is well known, and was used in the prior work on the subject, that if $u$ is a bounded subharmonic function, then the measure $\Delta u(z) \ln \frac1{|z|}\,dA(z)$ is Carleson. The new idea here is that one can get similar result for \emph{unbounded} subharmonic $u$ by multiplying the measure by an appropriate correcting factor, see Lemma \ref{l2.1} in the paper. While Lemma \ref{l2.1} itself is not needed for the proof, the idea of introducing correcting factors works for the embeddings on Hermitian vector bundles which we use in the proof, see Lemma \ref{l2.4}. 

The second idea is to use a more geometric approach to the problem (see Lemma \ref{l1.1} above), motivated by a surprising lemma by N.~Nikolski connecting solvablity of the corona problem with the existence of a bounded analytic projection, see Lemma 0.1 in \cite{TrWick}. This approach had allowed us to estimate \emph{only} embeddings involving the curvature term 
$$
\Delta \ln(\nm f(z)\nm^2)= \frac{\nm f\nm^2 \nm f'\nm^2 -|(f,f')|^2}{\nm f\nm^4 }, 
$$
while in the prior approaches, based on  modifications of T.~Wolff's proof of the Corona Theorem, the embeddings involving $\nm f'\nm^2/\nm f\nm^2$ were also required. Such embeddings apparently do not admit as good estimates as the ones with the curvature, which accounts for the extra $1/2$ in the exponent at $\ln s^{-1}$.  

\subsection{A conjecture and an open problem}
The statement of Theorem \ref{t0.1} looks like it is a final answer, so we conjecture that this theorem is sharp. Namely, let $\f\ge0$ be a bounded function on $[0,1]$ such that $\f(s)/s^2$ is an non-decreasing function. Any such function can be represented as $\f(s) =s^2 \psi(\ln s^{-2})$ where $\psi\ge0$ is a bounded  non-increasing function on $R_+$. 

\begin{conj}
If $\f(s) = s^2 \psi(\ln s^{-2})$ where $\psi\ge0$ is a bounded non-increasing function on $R_+$ such that 
$$
\int_0^\infty \psi(x) dx = \infty. 
$$
Then the condition 
$$
|\tau(z)|\le \f(\nm f(z)\nm ), \qquad \forall z\in \D
$$
does not imply that the equation $gf=\tau$ has a solution $g\in H^\infty_{E^*}$. 
\end{conj}

An interesting open problem would be to find a necessary and sufficient condition for the problem of ideals (the solvability of the equation $gf=\tau$). While it is pretty obvious that a size condition of the form $|\tau(z)|\le \f(\nm f(z)\nm)$ cannot be necessary and sufficient, the author still hopes that it is possible to find a necessary and sufficient condition in terms of, say, Carleson measures.  For example, that some measure constructed from the curvature  $\Delta \ln(\nm f(z)\nm^2) =\frac{\nm f\nm^2 \nm f'\nm^2 -|(f,f')|^2}{\nm f\nm^4 }$ and $\tau$ should be Carleson.  

The results of such type were obtained recently in \cite{TrWick}, where a necessary and sufficient condition for the solvability of the operator corona problem (left invertibility in $H^\infty$ of an operator-valued function $F$) were given. These condition were exactly the estimates on the curvature, namely the conditions that $\nm \p\Pi(z) \nm \le C (1-|z|)^{-1}$ and that the measure $\nm \p\Pi (z)\nm^2 \ln\frac1{|z|}\,dA(z)$ is Carleson.

\def\cprime{$'$}
 \def\lfhook#1{\setbox0=\hbox{#1}{\ooalign{\hidewidth
  \lower1.5ex\hbox{'}\hidewidth\crcr\unhbox0}}}
\providecommand{\bysame}{\leavevmode\hbox to3em{\hrulefill}\thinspace}
\providecommand{\MR}{\relax\ifhmode\unskip\space\fi MR }
\providecommand{\MRhref}[2]{%
  \href{http://www.ams.org/mathscinet-getitem?mr=#1}{#2}
}
\providecommand{\href}[2]{#2}
\providecommand{\bysame}{\leavevmode\hbox to3em{\hrulefill}\thinspace}


\begin{thebibliography}{10}

\bibitem{Andersson}M.~Andersson, {\em The Corona Theorem for Matrices}, 
{Math.~Z.},
{\bf 201} (1989), 121--130.

%

\bibitem{Bernd-db}
Bo~Berndtsson, \emph{{$\overline\partial\sb b$} and {C}arleson type
  inequalities}, Complex analysis, II (College Park, Md., 1985--86), Lecture
  Notes in Math., vol. 1276, Springer, Berlin, 1987, pp.~42--54.

\bibitem{Carl-corona1}
L.~Carleson, \emph{{Interpolations by bounded analytic functions and the
  Corona problem}}, Ann. of Math. (2) \textbf{76} (1962), 547--559.


\bibitem{Cegrell1990}
U.~Cegrell, \emph{A generalization of the corona theorem in the unit disc},
  Math. Z. \textbf{203} (1990), no.~2, 255--261. 

\bibitem{Cegrell1994}
\bysame, \emph{Generalisations of the corona theorem in the unit disc}, Proc.
  Roy. Irish Acad. Sect. A \textbf{94} (1994), no.~1, 25--30.
  

\bibitem{Gar-BAF}
J.~B.~Garnett, \emph{Bounded analytic functions}, Pure and Applied
  Mathematics, vol.~96, Academic Press Inc. [Harcourt Brace Jovanovich
  Publishers], New York, 1981.

\bibitem{Helson-lect}
H.~Helson, \emph{{Lectures on invariant subspaces}}, Academic Press, New
  York, 1964.

\bibitem{Lin-ideals-1993}
Kai-Ching Lin, \emph{On the constants in the corona theorem and the ideals of
  {$H\sp \infty$}}, Houston J. Math. \textbf{19} (1993), no.~1, 97--106.
  
\bibitem{Nagy-Foias-Book-1970}
B.~Sz.-Nagy and C.~Foia{\lfhook{s}}, \emph{Harmonic analysis of
  operators on {H}ilbert space}, Translated from the French and revised,
  North-Holland Publishing Co., Amsterdam, 1970. 
  
  
\bibitem{Nik-book-v1}
N.~K.~Nikolski, \emph{Operators, functions, and systems: an easy reading.
  {V}ol. 1: Hardy, Hankel, and Toeplitz}, Mathematical Surveys and
Monographs,
  vol.~92, American Mathematical Society, Providence, RI, 2002, Translated from
  the French by Andreas Hartmann.


\bibitem{Nik-shift}
\bysame, 
\emph{{Treatise on the shift operator}},
  Grundlehren der Mathematischen Wissenschaften [Fundamental Principles of
  Mathematical Sciences], vol. 273, Springer-Verlag, Berlin, 1986, Spectral
  function theory, With an appendix by S. V. Hru\v s\v cev [S. V.
  Khrushch{\"e}v] and V. V. Peller, Translated from the Russian by Jaak Peetre.


\bibitem{Pau2005}
J.~Pau, \emph{On a generalized corona problem on the unit disc}, Proc. Amer.
  Math. Soc. \textbf{133} (2005), no.~1, 167--174 (electronic). 
  
  

\bibitem{Rosenblum-corona}
M.~Rosenblum, \emph{A corona theorem for countably many functions},
  Integral Equations Operator Theory \textbf{3} (1980), no.~1, 125--137.


\bibitem{Rao67}
K.~V. Rajeswara~Rao, \emph{On a generalized corona problem}, J. Analyse Math.
  \textbf{18} (1967), 277--278.
  
  
\bibitem{Tol}
V.~A.~Tolokonnikov, \emph{{Estimates in the Carleson corona theorem,
ideals of
  the algebra $H\sp{\infty }$, a problem of Sz.-Nagy}}, Zap. Nauchn. Sem.
  Leningrad. Otdel. Mat. Inst. Steklov. (LOMI) \textbf{113} (1981), 178--198,
  267, Investigations on linear operators and the theory of functions, XI.
  
 
\bibitem{Treil-Ideals-2002}
S.~Treil, \emph{Estimates in the corona theorem and ideals of {$H\sp \infty$}:
  a problem of {T}.\ {W}olff}, J. Anal. Math. \textbf{87} (2002), 481--495,
  Dedicated to the memory of Thomas H.\ Wolff. 
  
\bibitem{TrlVol} S.~Treil, A.~Volberg, \emph{A Fixed Point Approach to Nehari's Problem and its Applications}, {Oper. Theory Adv. Appl., 71}, (1994), 165-186.

\bibitem{TrWick} S.~Treil, B.~Wick, \emph{Analytic projections, Corona Problem and geometry of holomorphic vector bundles}, Preprint, {\tt arXive:math.CA/0702756}.

\bibitem{Trent-Corona}
T.~Trent, \emph{{A new estimate for the vector valued corona
problem}},
  J. Funct. Anal. \textbf{189} (2002), no.~1, 267--282.
 
\bibitem{Uchiyama1981}
A.~Uchiyama, \emph{Corona theorems for countably many functions and
  estimates for their solution}, Preprint, USLA, 1981.
  
  


\end{thebibliography}
\end{document}